\documentclass{amsart}

\usepackage{amssymb,amsmath}
\usepackage{amsfonts,amsthm}
\newtheorem*{thm}{Theorem}

\newtheorem*{lem3.1}{\sc Lemma 3.1}
\newtheorem*{lem3.2}{\sc Lemma 3.2}
\newtheorem*{lem3.3}{\sc Lemma 3.3}
\newtheorem*{lem3.4}{\sc Lemma 3.4}
\newtheorem*{lem3.5}{\sc Lemma 3.5}
\newtheorem*{lem3.6}{\sc Lemma 3.6}
\newtheorem*{lem3.7}{\sc Lemma 3.7}
\newtheorem*{lem3.8}{\sc Lemma 3.8}
\newtheorem*{lem3.9}{\sc Lemma 3.9}

\begin{document}
\title{Smooth numbers in short intervals}
\author{K. Soundararajan} 
\thanks{The author is supported in part by a grant from the National Science Foundation}
\begin{abstract} 
Assume the Riemann Hypothesis. For every $\epsilon >0$ we 
show that there is a constant $C(\epsilon)$ such 
that for all large $x$, the interval $[x,x+C(\epsilon)\sqrt{x}]$ 
contains an integer all of whose prime factors are less 
than $x^{\epsilon}$. 
\end{abstract}
\maketitle
%\footnote{{\it AMS Subject Classification :} primary 11 M06.}   
%\footnote{{\it Key words :} Smooth numbers, Riemann zeta-function}
 
A natural number $n$ is called $y$-smooth if all its 
prime factors are below $y$.  We let ${\mathcal S}(y)$ denote 
the set of $y$-smooth numbers, and let $\Psi(x,y)$ denote the 
number of such integers below $x$.  If we write $y=x^{\frac 1u}$ 
then it is known that $\Psi(x,y) \sim \rho(u) x$ where $\rho(u)$ denotes 
the Dickman function defined by $\rho(u) =1$ for $0\le u\le 1$ and 
for $u\ge 1$ is defined as the unique continuous solution to 
the differential-difference equation $u{\rho}^{\prime}(u)=-\rho(u-1)$.  
This asymptotic formula was published first by Dickman \cite{D} for 
fixed values of $u$ and as $x\to \infty$; recently Soundararajan \cite{S} has 
pointed out that such an asymptotic formula may be found in Ramanujan's unpublished 
papers.  Later work has 
established asymptotic formulae for $\Psi(x,y)$ uniformly for $u$ in a 
wide range; see for example  the surveys \cite {G} and \cite{HT}. 

In this note we are concerned with the existence of smooth numbers in 
short intervals.   For a wide range of the variables $x$, $y$, and $z$,
it is expected that
\begin{equation*}
\label{E1}
\Psi (x+z, y) - \Psi (x, y) \asymp \frac {z}{x} \Psi (x,y).
\end{equation*}
It is also of interest to establish the existence of smooth numbers in 
such short intervals, even if one is not able to exhibit a positive proportion of 
such numbers.  One motivation for this problem is the analysis of Lenstra's 
elliptic curve factorization algorithm \cite{L} (and see also \cite{P}) where one wishes to find 
integers in $[x,x+4\sqrt{x}]$ which are $\exp(\sqrt{\log x\log \log x})$ smooth.  
  
  Regarding this problem, an important advance was made by Balog \cite{B}
  who showed that for any fixed $\epsilon >0$ and $x$ large, the interval
  $[x,x+x^{\frac 12+\epsilon}]$ contains many $x^{\epsilon}$-smooth integers.   
  Harman \cite{Ha1} has obtained a strengthening of this result, allowing $\epsilon$ to 
  be a function of $x$.  The problem for intervals of length of size $\sqrt{x}$ 
 has proved resistant, and here Harman \cite{Ha2} has established that the interval $[x,x+\sqrt{x}]$ 
 contains integers that are $x^{\frac 1{4\sqrt{e}}}$-smooth.  If one expands the 
 interval a little to consider $[x, x+C\sqrt{x}]$ for some constant $C>0$, then 
 recently Matom{\" a}ki \cite{M}, advancing an approach of Croot \cite{C},  has shown that such intervals (for 
 a suitably large value of $C$) contain $x^{\frac 1{5\sqrt{e}}+\epsilon}$-smooth numbers.  
 
 One may wonder if the Riemann Hypothesis is of use in this problem.   Assuming RH, 
 Xuan \cite{X} has shown that intervals $[x,x+x^{\frac 12}(\log x)^{1+\epsilon}]$ contain 
 $x^{\epsilon}$-smooth integers.   Recently Ganguly and Pal \cite{GP} have noted that if,
 in addition to RH, one assumes a strong conjectural estimate for $\pi S(t)$ (which is  
 the argument of $\zeta(1/2+it)$), then intervals $[x,x+x^{\frac 12} (\log x)^{\frac 12+\epsilon}]$ contain 
 $x^{\epsilon}$-smooth numbers.  
 We improve upon Xuan's work by establishing the following theorem, 
 which unfortunately is still not strong enough to be applicable to the analysis of Lenstra's algorithm.

\begin{thm}
Assume the Riemann Hypothesis.  Let $x$ be large and suppose that $x\ge y\ge \exp(5\sqrt{\log x\log \log x})$, and write $y=x^{\frac 1u}$.   There is an absolute constant $B$ such that 
with $z = Bu\sqrt{x}/\rho(u/2)$ we have 
$$ 
\Psi(x+z, y) - \Psi(x,y) \gg_{\epsilon} zx^{-\epsilon}. 
$$ 
\end{thm}

By using estimates for divisor functions in short intervals we can obtain a better 
lower bound for the number of smooth integers in short intervals, but our methods 
would not give a positive proportion.  Our Theorem sheds no light 
on $y$-smooth integers with $y$ smaller than $\exp(\sqrt{\log x\log \log x})$, 
and it would be interesting to devise alternative approaches in this regime.  
Our Theorem is also likely to be very far from the truth about smooth numbers 
in short intervals.  For example, one would expect that for every $\epsilon >0$ there 
exists a constant $C(\epsilon)$ such that every interval $[x,x+C(\epsilon)\log x]$ 
contains an $x^{\epsilon}$-smooth number.  This would be analogous to 
Cram{\' e}r's conjecture on the distribution of prime numbers, and note that 
a similarly large gulf exists between what can be established about primes on RH and 
the expected truth.  As with primes, one can say more about the existence of 
smooth numbers in almost all short intervals; this problem has been considered by 
Hafner \cite{Haf} but his work remains unpublished.  

We now turn to the proof of our Theorem.  Let the parameters $x$, $y$, $z$ and $u$ be as in the Theorem, and define $\delta$ by $xe^{2\delta} = x+z$.  Let 
$$ 
M(s) = \sum_{{\sqrt{x}y^{-1/3} \le n\le \sqrt{x}y^{-1/4}} \atop {n\in {\mathcal S}(y)} } 
\frac 1{n^s}. 
$$ 
Our proof of the Theorem is based upon considering 
$$ 
I = \frac{1}{2\pi i} \int_{c-i\infty}^{c+i\infty} -\frac{\zeta^{\prime}}{\zeta}(s) M(s)^2 x^s 
\frac{(e^{\delta s}-1)^2}{s^2 } ds  
$$ 
where $c=1+\frac {1}{\log x}$.

By shifting contours to the left if $\xi >1$ and to the right if $\xi \le 1$ we may 
see that 
$$ 
\frac{1}{2\pi i} \int_{c-i\infty}^{c+i\infty} \xi^s \frac{ds}{s^2} = 
\begin{cases}
\log \xi &\text{if } \xi \ge 1\\
0 &\text{if  } 0 <\xi \le 1.\\
\end{cases} 
$$ 
Therefore 
$$ 
\frac{1}{2\pi i} \int_{c-i\infty}^{c+i\infty} 
\xi^s \frac{(e^{\delta s}-1)^2}{s^2} ds 
= 
\begin{cases}
\min( \log(e^{2\delta}\xi), \log (1/\xi))  &\text{if } e^{-2\delta} \le \xi \le 1 \\ 
0 &\text{otherwise}. \\
\end{cases}
$$ 
Hence
$$
I= \sum_{x \le n \le x e^{2\delta}}  \sum_{{{n= rm_1 m_2} \atop {m_1, m_2 \in {\mathcal S}(y)} }
\atop {\sqrt{x}y^{-1/3 }\le m_1, m_2 \le \sqrt{x}y^{-1/4} } } \Lambda(r) 
\min \Big( \log \frac{e^{2\delta}x}{n}, \log \frac{n}{x}\Big).  
$$
Note that $r=n/(m_1 m_2)$ is at most $y$, and so the integers $n$ counted in the 
RHS above are all $y$-smooth.  Moreover the number of ways of writing $n$ as $rm_1m_2$ is 
at most $d_3(n) \ll x^{\epsilon}$.  Therefore we conclude that 
\begin{equation}
\label{UB}
I \ll \delta x^{\epsilon}  ( \Psi(x+z,y)-\Psi(x,y)).
\end{equation}

We shall now derive a lower bound for $I$ which will prove the Theorem.  We 
move the line of integration in the definition of $I$ to the line Re$(s)=-\frac 12$.  
We encounter poles at $s=1$ and at the non-trivial zeros $\rho=\frac 12+i\gamma$ 
of $\zeta(s)$.  Thus we find that $I$ equals 
\begin{equation}
\label{LB1}
x (e^{\delta}-1)^2 M(1)^2 - \sum_{\rho} M(\rho)^2 x^{\rho} 
\Big( \frac{e^{\delta \rho}-1}{\rho}\Big)^2 - \frac{1}{2\pi i } 
\int_{-\frac 12-i\infty}^{-\frac 12+i\infty} \frac{\zeta^{\prime}}{\zeta}(s) M(s)^2 
x^s \frac{(e^{\delta s}-1)^2}{s^2} ds. 
\end{equation}

 Using the functional equation for 
 $\zeta(s)$ and Stirling's formula, we obtain that $|\frac{\zeta^{\prime}}{\zeta}(-\tfrac 12+it)| 
 \ll \log (2+|t|)$.  Since $(e^{\delta s}-1)^2/s^2 \ll \min( \delta^2, 1/|s|^2)$ for all 
 complex $s$ with $-\frac 12 \le \text{Re }(s)\le 2$, we find that the integral 
 appearing in \eqref{LB1} is bounded by 
 $$ 
 \ll x^{-\frac 12} \int_{-\infty}^{\infty} |M(-\tfrac 12 +it)|^2 \log (2+|t|) 
 \min \Big( \delta^2, \frac{1}{1/4+t^2}\Big) dt. 
 $$ 
 We now split the interval $(-\infty, \infty)$ into the sets 
 ${\mathcal I}_0= \{ |t|\le 1/\delta\}$, and ${\mathcal I}_j 
 = \{ 2^{j-1} /\delta \le |t| \le 2^{j}/\delta\}$ for $j\in {\Bbb N}$.
By appealing to a standard mean-value theorem for Dirichlet 
polynomials (see for example Theorem 9.1 of \cite{IK}) we find that the 
contribution from $t \in {\mathcal I}_0$ is
\begin{eqnarray*}
&\ll& x^{-\frac 12} \delta^2 \log (1/\delta) 
\int_{{\mathcal I}_0} |M(-\tfrac 12+it)|^2 dt 
\\
&\ll& x^{-\frac 12} \delta^2 \log x \Big( \frac{\sqrt{x}}{y^{\frac 14}} + \frac{1}{\delta}\Big) 
\sum_{{n \in {\mathcal S}(y)} \atop {\sqrt{x}/y^{\frac 13} \le n\le \sqrt{x}/y^{\frac 14}}} n \\
&\ll&  \delta^2 \log x \frac{\sqrt{x}}{y^{\frac 12}} \Big( \frac{\sqrt{x}}{y^{\frac 14}} + \frac 1\delta\Big) M(1).\\
\end{eqnarray*}
Similarly we find that for $j\ge 1$ the contribution from the interval ${\mathcal I}_j$ is 
$$ 
\ll \frac{j \delta^2 \log x)}{2^j} \frac{\sqrt{x}}{\sqrt{y}}\Big( \frac{\sqrt{x}}{y^{\frac 14}} + \frac{1}{\delta}
\Big) M(1).
$$ 
 We conclude that the integral in \eqref{LB1} is bounded by 
\begin{equation}
\label{LB2}
\ll  \delta^2 \log x\frac{\sqrt{x}}{\sqrt{y}} \Big( \frac{\sqrt{x}}{y^{\frac 14}} + \frac 1\delta\Big) 
M(1).
 \end{equation}
 
 Now we turn to the sum over zeros in \eqref{LB1}.  This sum is 
 bounded by 
 $$ 
 \ll x^{\frac 12} \sum_{\gamma} |M(\tfrac 12+i\gamma)|^2 \min \Big( \delta^2, 
 \frac{1}{1/4+\gamma^2}\Big).
 $$ 
 To estimate this, we decompose the sum into cases depending on whether
 $\gamma \in {\mathcal I}_j$ with ${\mathcal I}_j$ as earlier.   We shall prove that 
 the contribution from the zeros in ${\mathcal I}_j$ for any $j\ge 0$ is 
 \begin{equation}
 \label{LB3}
 \ll  \sqrt{x} \frac{j+1}{2^j} \delta^2 \log (1/\delta) \Big( \frac{\sqrt{x}}{y^{\frac14}} + \frac{1}{\delta}\Big) M(1).
 \end{equation} 
 Summing over all $j$, it then follows that the sum over zeros in \eqref{LB1} is
 \begin{equation} 
 \label{LB4} 
 \ll \sqrt{x}\delta^2 \log (1/\delta) \Big( \frac{\sqrt{x}}{y^{\frac14}} + \frac{1}{\delta}\Big) M(1). 
 \end{equation}

 Note that the contribution in \eqref{LB2} is dominated by that in \eqref{LB4}.  Thus, combining \eqref{LB1}, \eqref{LB2} and \eqref{LB4} we conclude that 
 $$ 
 I  \ge x \delta^2 M(1)^2 -A\sqrt{x} \delta^2 \log x \Big(\frac{\sqrt{x}}{y^{\frac 14}} 
 + \frac 1\delta\Big) M(1),
 $$ 
 for an appropriate absolute constant $A$.  Since $\rho(u/2) = u^{-u/2(1+o(1)}$, in our range of $\delta$ and $y$ we have that $1/\delta \ge \sqrt{x}/y^{\frac 14}$.  Moreover, using the 
 asymptotic formula for smooth numbers, 
 we see that $M(1) \ge \rho(u/2) (\log y)/24$.  Choosing $B$ suitably 
 in terms of $A$, from the above remarks we find 
 that  $I\ge x\delta^2 M(1)^2/2$, and by \eqref{UB}
  the Theorem follows.

 It remains lastly to justify the bound \eqref{LB3}.  We treat the case $j=0$, the other cases being 
 similar.   The proof is entirely standard and we sketch the details quickly; indeed we could 
 obtain asymptotic formulae for such sums but we do not need this.  Let $\xi(s) = s(s-1)\pi^{-s/2} \Gamma(s/2)\zeta(s)$ denote Riemann's $\xi$-function 
 which is entire, and whose zeros are the non-trivial zeros of the $\zeta(s)$.  We consider, 
 with $c=1+1/\log x$, 
 $$ 
 J := \frac{1}{2\pi i} \int_{(c)} \frac{\xi^{\prime}}{\xi}(s) M(s) M(1-s) \Big( \frac{e^{\delta(s-1/2)}-e^{-\delta(s-1/2)}}{(s-1/2)}\Big)^2 ds.
 $$ 
 We now move the line of integration to Re$(s)=1-c$.  There are poles at the non-trivial 
 zeros of $\zeta(s)$ and these contribute 
 $$ 
 \sum_{\gamma} |M(1/2+i\gamma)|^2 \Big(\frac{2\sin(\delta \gamma)}{\gamma}\Big)^2.
 $$ 
 To handle the integral on the line Re$(s)=1-c$ we use the functional equation 
 $\xi^{\prime}/\xi(s) = -\xi^{\prime}/\xi(1-s)$ and then make a change of variable $w=1-s$.  
 In this manner we recognize the integral on Re$(s)=1-c$ as being $-J$.   
 Thus we find that 
 $$ 
 2J  = \sum_{\gamma} |M(1/2+i\gamma)|^2 \Big(\frac{2\sin(\delta \gamma)}{\gamma}\Big)^2 
 \gg \delta^2 \sum_{|\gamma|\le 1/\delta} |M(1/2+i\gamma)|^2. 
 $$ 
We may therefore focus on bounding $J$.  

Note that 
$$ 
\frac{\xi^{\prime}}{\xi}(s) = \Big( \frac{1}{s}+\frac{1}{s-1} -\log \sqrt{\pi} +\frac 12\frac{\Gamma^{\prime}}{\Gamma}(s/2)\Big) + \frac{\zeta^{\prime}}{\zeta}(s),
$$ 
and accordingly write $J=J_1 + J_2$.  To estimate $J_2$, we expand the Dirichlet 
series for $\zeta^{\prime}/\zeta(s)$, $M(s)$ and $M(1-s)$ and exchange the summations 
and integration.  Thus 
$$ 
J_2 = -\sum_{n=1}^{\infty} \Lambda(n) \sum_{{m_1, m_2 \in {\mathcal S}(y)}  
\atop {\sqrt{x}/y^{1/3} \le m_1, m_2 \le \sqrt{x}/y^{\frac 14}} } \frac{1}{m_2}
\frac{1}{2\pi i} \int_{(c)} \Big(\frac{m_2}{nm_1}\Big)^s \Big(\frac{e^{\delta(s-1/2)}-e^{-\delta(s-1/2)}}{(s-1/2)} \Big)^2 ds. 
$$ 
By shifting contours appropriately, we find that 
$$ 
\frac{1}{2\pi i} \int_{(c)} \xi^{s} \Big(\frac{e^{\delta(s-1/2)} -e^{-\delta(s-1/2)}}{(s-1/2)}\Big)^2 ds
= 
\begin{cases} 
\sqrt{\xi} (2\delta -|\log \xi|) &\text{if  } e^{-2\delta}\le \xi \le e^{2\delta}\\ 
0 &\text{otherwise}.\\
 \end{cases}
 $$ 
 Since $m_1$ and $m_2$ are below $\sqrt{x}/y^{\frac 14}$, if $nm_1 \neq m_2$ then 
 $m_2/nm_1$ lies outside the interval $(e^{-2\delta},e^{2\delta})$.  Thus 
 $$ 
 J_2 = - 2\delta \sum_{{m_2\in {\mathcal S}(y)} \atop {\sqrt{x}/y^{1/3} 
 \le m_2\le \sqrt{x}/y^{1/4}} } \frac{1}{m_2} \sum_{{nm_1=m_2} \atop {\sqrt{x}/y^{1/3} \le m_1 
 \le \sqrt{x}/y^{1/4}}}\Lambda(n) <0. 
 $$  
 Thus $J\le J_1$, and we are reduced to estimating $J_1$.  
 
 To estimate $J_1$ we move the line of integration to the line Re$(s)=1/2$.  We 
 encounter a pole at $s=1$ whose residue is $4M(0)M(1)(e^{\delta/2}-e^{-\delta/2})^2$.  
 Using Stirling's formula we find that the remaining integral on the Re$(s)=1/2$ line 
 is 
 $$ 
 \ll \int_{-\infty}^{\infty} |M(1/2+it)|^2 \log (2+|t|) \Big(\frac{\sin(\delta t)}{t}\Big)^2 dt. 
 $$ 
 Splitting this integral into the intervals ${\mathcal  I}_j$ as above, and appealing to the mean value theorem for Dirichlet polynomials we conclude that this quantity is 
 $$ 
 \ll \delta^2 \log x \Big(\frac{\sqrt{x}}{y^{\frac 14}}+\frac 1\delta\Big) M(1). 
 $$
 Since $M(0) \le \sqrt{x}/y^{\frac 14}$ we conclude that 
 $$ 
J\le  J_1 \ll \delta^2 \log x \Big( \frac{\sqrt{x}}{y^{\frac 14}} +\frac 1\delta\Big) M(1).
 $$ 
 This proves \eqref{LB4} for the region ${\mathcal I}_0$, and as noted before the 
 other cases follow similarly.  Our proof of the Theorem is now complete.
  
\vskip .1 in 

\noindent {\bf Acknowledgments.}  This note was inspired by an unpublished manuscript 
of Ganguly, Pal, and Sankaranarayanan. %where they strengthened Xuan's work 
%and established that intervals $[x,x+C(\epsilon) \sqrt{x}\log x/\log \log x]$ 
%contained $x^{\epsilon}$-smooth numbers.  
I offered them joint authorship of this 
note, but they declined and preferred that I write it up separately.  I thank them for their 
comments and valuable correspondence.   I also thank Carl Pomerance for drawing my attention 
to Hafner's work \cite{Haf}.

\vskip .25 in

 \noindent  Department of Mathematics, \\
 Stanford University, \\
 450 Serra Mall, Building 380, \\
 Stanford,  CA 94305-2125, USA.\\
{ksound@math.stanford.edu} \\

\end{document}